\documentclass{ita}

\makeatletter
\renewcommand\normalsize{\@xsetfontsize\normalsize 6%
  \@adjustvertspacing \let\@listi\@listI 
  \abovedisplayskip 9pt \@plus2pt \@minus2pt 
  \belowdisplayskip \abovedisplayskip}
\def\@NMJRNL@F{{\relax}}%
\def\@NMJRNL@E{{\relax}}%
\def\@JRNL@X{{\relax}}%
\def\ps@firstpage{\ps@plain
  \def\@oddfoot{\hfill{\scriptsize \copyright\ EDP Sciences 2012}}%
  \let\@evenfoot\@oddfoot\def\@oddhead{\null\hss}
  \let\@evenhead\@oddhead}
\numberwithin{equation}{section}
\editor{G. Richomme}
\def\@date{Recieved November 2, 2010. Accepted February 3, 2012.}
\def\@MSSG@RNNGTTL{\MakeUppercase{Morphisms preserving the set of 3-interval exchange words}}
\def\@MSSG@RNNGTHR{T. Hejda}
\makeatother

\usepackage{amssymb}

\begin{document}

\addtocounter{page}{106}

\noindent\begingroup\small
\begin{tabular}{@{}p{\textwidth}@{}}
RAIRO-Theor. Inf. Appl. 46 (2012) 107--122
\hfill
Available online at:
\\
DOI: 10.1051/ita/2012009
\hfill
www.rairo-ita.org
\end{tabular}
\endgroup
\vspace{1.3cm}\leavevmode

\title{Morphisms preserving the set of words coding three interval exchange}
\runningtitle{Morphisms preserving the set of 3-interval exchange words}
\runningauthors{Tom\'a\v s Hejda}

\thanks{
We acknowledge financial support by the
     Czech Science Foundation grant \emph{201/09/0584}
and by the grants
     \emph{MSM6840770039} and \emph{LC06002} of the Ministry of Education, Youth, and Sports of the Czech Republic.
We also thank the
     CTU student grant \emph{SGS10/085/OHK4/1T/14} and \emph{SGS11/162/OHK4/3T/14}.
}
\thanks{
We would like to thank the organizers of the conference
 \emph{13$^{\text{i\`emes}}$ Journ\'ees Montoises d'Informatique Th\'eorique}
 for a financial support of the author's stay at the conference.
}

\author{Tom\'a\v s Hejda}
\address{
Department of Mathematics FNSPE,
Czech Technical University in Prague,
Trojanova 13, 120 00 Prague, Czech Republic;
e-mail: \texttt{tohecz@gmail.com}
}



\begin{abstract}
Any amicable pair $\varphi$, $\psi$ of Sturmian morphisms enables a construction of a ternary morphism
 $\eta$ which preserves the set of infinite words coding 3-interval exchange.
We determine the number of amicable pairs with the same incidence matrix in $\mathrm{SL}^\pm(2,\mathbb{N})$
 and we study incidence matrices associated with the corresponding ternary morphisms $\eta$.
\end{abstract}

\maketitle

\section{Introduction}

\emph{Sturmian words} are well-described objects in combinatorics on words.
They can be defined in several equivalent ways \cite{Ber96}, e.g.\ as words coding a two-interval exchange
 transformation with irrational ratio of lengths of the intervals.
Morphisms preserving the set of Sturmian words are called \emph{Sturmian} and they form
 a monoid generated by three of its elements (see \cite{BS94,Loth}).
Let us denote this monoid by $\mathcal{M}_{\mathrm{Sturm}}$.

In this paper, we consider morphisms preserving the set of words coding a three-interval
 exchange transformation with permutation $(3,2,1)$, the so-called \emph{3iet words}.
We call these morphisms \emph{3iet-preserving}.
Monoid of these morphisms, denoted by $\mathcal{M}_{\mathrm{3iet}}$, is not fully described.
It is shown (see \cite{Hak}) that the monoid $\mathcal{M}_{\mathrm{3iet}}$ is not finitely generated.
Recently, in \cite{AMP-Mor}, pairs of amicable Sturmian morphisms were defined.
The authors used this notion to describe morphisms that have as a fixed point
 a non-degenerate 3iet word, i.e. word with complexity $\mathcal{C}(n)=2n+1$.
Using the operation of ``ternarization'', we can assign a morphism $\eta=\operatorname{ter}(\varphi,\psi)$ over
a ternary alphabet to a pair of amicable Sturmian morphisms.
We show that such $\eta$ is a 3iet-preserving morphism.
Moreover, we show that the set
\[
	\mathcal{M}_{\mathrm{ter}}=\bigl\{\operatorname{ter}(\varphi,\psi)\big|\varphi,\psi \text{ amicable morphisms}\bigr\}
\]
 is a monoid, but it does not cover the whole monoid $\mathcal{M}_{\mathrm{3iet}}$.

We also study the incidence matrices of morphisms $\eta\in\mathcal{M}_{\mathrm{ter}}$.
From the definition of amicable Sturmian morphisms $\varphi,\psi$ we can derive that $\varphi$ and $\psi$ have
 the same incidence matrix $\mathbf A\in\mathbb{N}^{2\times 2}$, where $\det\mathbf A=\pm1$.
As shown in~\cite{Seebold}, for every matrix $\mathbf A=\left(\begin{smallmatrix} p_0 & q_0 \\ p_1 & q_1 \end{smallmatrix}\right)$
 with $\det\mathbf A=\pm1$, there exist $p_0+p_1+q_0+q_1-1$ Sturmian morphisms.
We will show the following theorem concerning the number of pairs of amicable Sturmian morphisms with
 a given matrix.

\begin{thrm}\label{thm:AmN}
	Let $\mathbf A=\left(\begin{smallmatrix} p_0 & q_0 \\ p_1 & q_1 \end{smallmatrix}\right)\in\mathbb{N}^{2\times2}$ be a matrix with $\det\mathbf A=\pm1$.
	Then there exist exactly
	\begin{equation}\label{eq:AmN}
		m\bigl(\left\|\mathbf A\right\|-1\bigr)+\frac m2\bigl(\det\mathbf A-m\bigr)
	\end{equation}
	 pairs of amicable Sturmian morphisms with incidence matrix $\mathbf A$,
	 where $m=\min\{p_0+p_1, q_0+q_1\}$ and $\left\|\mathbf A\right\|=p_0+p_1+q_0+q_1$.
\end{thrm}

Moreover, for a given matrix $\mathbf A$, we will describe all matrices $\mathbf B\in\mathbb{N}^{3\times3}$ such that
 $\mathbf B$ is an incidence matrix of $\eta=\operatorname{ter}(\varphi,\psi)$ for amicable Sturmian morphisms
 $\varphi, \psi$ with incidence matrix $\mathbf A$.

\section{Preliminaries}

\subsection{Words over finite alphabet}
Besides the infinite words, we consider \emph{finite words} over the alphabet $\mathbb{A}$.
We write $w=w_0w_1\cdots w_{n-1}$, where $w_i\in\mathbb{A}$ for all $i\in\mathbb{N}$, $i<n$.
We denote by $\left|w\right|$ the length $n$ of the finite word $w$.
We denote by $\left|w\right|_a$ the number of occurrences of a letter $a\in\mathbb{A}$ in the word $w$.
The set of all finite words on the alphabet $\mathbb{A}$ including the empty word is denoted by $\mathbb{A}^*$.
The set $\mathbb{A}^*$ with the operation of concatenation is a monoid.
On the set $\mathbb{A}^*$ we define a relation of \emph{conjugation}:
 $w\sim w'$, if there exists $v\in\mathbb{A}^*$ such that $wv=vw'$.
A \emph{morphism} from $\mathbb{A}^*$ to $\mathcal B^*$ is a mapping $\varphi:\mathbb{A}^*\rightarrow\mathcal B^*$ such that
 $\varphi(vw)=\varphi(v)\varphi(w)$ for all $v,w\in\mathbb{A}^*$.
It is clear that a morphism is well defined by images of letters $\varphi(a)$ for all $a\in\mathbb{A}$.
If $\mathbb{A}=\mathcal B$, then $\varphi$ is called a \emph{morphism over} $\mathbb{A}$.

The set of \emph{infinite words} over the alphabet $\mathbb{A}$ is denoted by $\mathbb{A}^\mathbb{N}$.
The action of a morphism can be naturally extended to an infinite word $(u_i)_{i\in\mathbb{N}}$ putting
 $\varphi(u)=\varphi(u_0)\varphi(u_1)\varphi(u_2)\cdots$.
If an infinite word $u\in\mathbb{A}^\mathbb{N}$ satisfies $\varphi(u)=u$, we call it a \emph{fixed point}
 of the morphism $\varphi$ over $\mathbb{A}$.

To a morphism $\varphi$ over $\mathbb{A}$ we assign an \emph{incidence matrix} $\mathbf M_\varphi$ defined by
 $(\mathbf M_\varphi)_{ab} = \left|\varphi(a)\right|_b$ for all $a,b\in\mathbb{A}$.
To a finite word $v\in\mathbb{A}^*$ we assign a \emph{Parikh vector} $\Psi(v)$ defined by
 $\Psi(v)_b = \left|v\right|_b$ for all $b\in\mathbb{A}$.

The \emph{language} of an infinite word $u$ is the set of all its factors.
Let us recall that a finite word $w\in\mathbb{A}^*$ is a \emph{factor} of $u=(u_i)_{i\in\mathbb{N}}$,
 if there exist indices $n,j\in\mathbb{N}$ such that $w=u_nu_{n+1}\cdots u_{n+j-1}$.
The language of an infinite word is denoted by $\mathcal L(u)$.

It is known that the language of neither Sturmian nor 3iet word depends on the point
 $x_0\in[0,1)$, the orbit of which the infinite word codes.
It depends only on slope $\varepsilon$ or parameters $\alpha,\beta$.

The \emph{(factor) complexity} of an infinite word $u$ is a mapping $\mathcal{C}_u:\mathbb{N}\rightarrow\mathbb{N}$,
 which returns the number of factors of $u$ of the length $n$, thus
 $\mathcal{C}_u(n)=\#\bigl\{w\in\mathcal L(u)\big|\left|w\right|=n\bigr\}$.
It is easy to see that a word $u$ is periodic if and only if there exists $n_0\in\mathbb{N}$
 such that $C_u(n_0)\leq n_0$.

\subsection{Interval exchange}
We consider Sturmian words, i.e.\ aperiodic words given by exchange of 2 intervals with
 permutation $(2,1)$, and words given by exchange of 3 intervals with permutation $(3,2,1)$.
Let us recall that general $r$-interval exchange transformations were introduced already in \cite{Katok67}.

The 2-interval exchange transformation $S$ is a mapping $S:[0,1)\rightarrow[0,1)$.
It is determined by its slope $\varepsilon\in[0,1]$ and is given by
\[
	Sx=\begin{cases}
	x+1-\varepsilon &\text{if }\; x\in[0,\varepsilon)
	\\
	x-\varepsilon &\text{if }\; x\in[\varepsilon,1)
	. \end{cases}
\]
The orbit of a point $x_0\in[0,1)$ with respect to the transformation $S$, i.e.\ the sequence
 $x_0, Sx_0, S^2x_0,\ldots$ can be coded by an infinite word $u=(u_i)_{i=0}^\infty$ on the
 binary alphabet $\{0,1\}$.
The infinite word is given by
\begin{equation}\label{eq:2iet-code}
	u_i=\begin{cases}
	0 &\text{if }\; S^ix_0\in[0,\varepsilon)
	,\\
	1 &\text{if }\; S^ix_0\in[\varepsilon,1)
	. \end{cases}
\end{equation}
It is a well-known fact that for an irrational $\varepsilon$, the word $u$ is Sturmian.
Using the same construction on the partition of the interval $(0,1]$ into
 $(0,\varepsilon]\cup(\varepsilon,1]$, we again obtain a Sturmian word.
On the other hand, every Sturmian word can be obtained by one of the above two constructions.
The set of Sturmian words will be denoted by $\mathcal{W}_{\mathrm{Sturm}}$.

In \cite{Loth} (the original results can be found in \cite{CoHe,HeMo}),
 the authors show that Sturmian words are the aperiodic words with minimal complexity,
 i.e. $\mathcal{C}_u(n)=n+1$ for all $u\in\mathcal{W}_{\mathrm{Sturm}}$ and $n\in\mathbb{N}$.
We can see that
\begin{equation}\label{eq:Sxi}
	S^ix_0=\{x_0-i\varepsilon\}
	\quad\text{for all}\quad
	x_0\in[0,1)
,\end{equation}
 where $\{x\}=x-\lfloor x\rfloor$ denotes the \emph{fractional part} of a number $x\in\mathbb R$.
Then $u_i=\lfloor x_0-i\varepsilon\rfloor-\lfloor x_0-(i+1)\varepsilon\rfloor$, which is exactly the formula how
 \cite{Loth} define mechanical words.
 
We will use another fact about the two-interval exchanges.
Let $\varphi\in\mathcal{M}_{\mathrm{Sturm}}$ be a Sturmian morphism.
Then the word $v=\varphi(a)$ for $a\in\{0,1\}$ codes two-interval exchange with the slope
 $\frac{\left|v\right|_0}{\left|v\right|}$.
We should see this from \cite[Lemma 2.1.15]{Loth}.
The word $a^k$ is a factor of some Sturmian word, hence the word $\varphi(a)^k$
 is balanced for any $k\in\mathbb{N}$, which means that the infinite word
 $u=\varphi(a)^\omega=\varphi(a)\varphi(a)\varphi(a)\cdots$ is balanced and periodic,
 thus it is rational mechanical.
In our terms, this means that it codes a rational 2-interval exchange;
 it is as well shown there that the slope of the transformation is exactly $\frac{\left|v\right|_0}{\left|v\right|}$.

\enlargethispage{3mm}

The 3-interval exchange transformation $T$ is determined by two parameters $\alpha,\beta\in(0,1)$
 satisfying $\alpha+\beta<1$.
Using parameters $\alpha$, $\beta$ and $\gamma=1-\alpha-\beta$ we partition the interval $[0,1)$
 into $I_A=[0,\alpha)$, $I_B=[\alpha,\alpha+\beta)$ and $I_C=[\alpha+\beta,1)$.
The mapping $T$ is given by
\[
	Tx=\begin{cases}
	x+\beta+\gamma &\text{if }\; x\in I_A
	,\\
	x-\alpha+\gamma &\text{if }\; x\in I_B
	,\\
	x-\alpha-\beta &\text{if }\; x\in I_C
	. \end{cases}
\]
The orbit of a point $x_0\in[0,1)$ with respect to the transformation $T$ is coded by a word
 $u=(u_i)_{i=0}^\infty$ over the ternary alphabet $\{A,B,C\}$:
\[
	u_i=X \quad\text{if}\quad T^i x_0\in I_X
. \]
Similarly to the case of 2-interval exchange transformation, we can define the exchange of 3 intervals using
the partition $(0,1]=(0,\alpha]\cup(\alpha,\alpha+\beta]\cup(\alpha+\beta,1]$.
If $\frac{1-\alpha}{1+\beta}$ is irrational, the infinite word $u$ is aperiodic, and we call it a
 \emph{3iet word}; the set of these words is denoted by $\mathcal{W}_{\mathrm{3iet}}$.
For combinatorial properties of 3iet words, see \cite{Ferenczi}.

Aperiodic words coding 3-interval exchange transformations, called here 3iet words, have the complexity
 $\mathcal{C}_u(n)\leq 2n+1$ for all $n\in\mathbb{N}$.
If a 3iet word $u\in\mathcal{W}_{\mathrm{3iet}}$ satisfies $\mathcal{C}_u(n)=2n+1$ for all $n\in\mathbb{N}$, we call it
 a \emph{non-degenerate} 3iet word;
 otherwise we call it a \emph{degenerate} 3iet word and it is a quasi-Sturmian word
 (see \cite{Cass}).

\subsection{Standard pairs and standard morphisms}
In \cite{Seebold}, the notion of standard pairs is introduced.
If we define two operators on pairs of words $L,R:\{0,1\}^*\times\{0,1\}^*\rightarrow\{0,1\}^*\times\{0,1\}^*$ as
\[
	L(x,y)=(x,xy)
	,\qquad
	R(x,y)=(yx,y)
,\]
we say that a pair $(x,y)$ is a \emph{standard pair}, if it can be obtained from the pair $(0,1)$ by applying
 the operators $L$ and $R$ finitely many times.
For every standard pair $(x,y)$ there exists a word $v\in\{0,1\}^*$ such that
\begin{equation}\label{eq:std-w01}
	xy=v01
\quad\text{and}\quad
	yx=v10
.\end{equation}

We say that a binary morphism $\varphi$ is \emph{standard},
 if there exists a standard pair $(x,y)$ such that
\[
	\begin{aligned}
		\varphi(0)&=x
		,\\
		\varphi(1)&=y
	,\end{aligned}
	\qquad\quad\text{or}\qquad\quad
	\begin{aligned}
		\varphi(0)&=y
		,\\
		\varphi(1)&=x
	.\end{aligned}
\]

The authors of \cite{Seebold} show the close connection between the standard morphisms
 and all the Sturmian morphisms:
\begin{enumerate}
 
\item
Every standard morphism is Sturmian.

\item
For every matrix $\mathbf A\in\mathbb{N}^{2\times2}$ with $\det\mathbf A=\pm1$,
 there exists exactly one standard morphism $\varphi$ with incidence matrix $\mathbf M_\varphi=\mathbf A$.

\item
Every Sturmian morphism $\psi\in\mathcal{M}_{\mathrm{Sturm}}$ is a right conjugate to some standard morphism $\varphi$.
Let us recall that a morphism $\psi$ over $\mathbb{A}$ is a \emph{right conjugate} to $\varphi$,
 if there exists a finite word $v\in\mathbb{A}^*$ such that
\[
	\varphi(a)v=v\psi(a)
	\quad\text{for all letters}\quad
	a\in\mathbb{A}
.\]

\end{enumerate}

\subsection{Amicable words and morphisms}
In the article \cite{ABMP}, authors show the close connection between 3iet and Sturmian words using
 morphisms $\sigma_{01},\sigma_{10}:\{A,B,C\}^*\rightarrow\{0,1\}^*$ given by
\begin{align*}
	\sigma_{01}(A)&=0    ,&   \sigma_{10}(A)&=0    ,\\
	\sigma_{01}(B)&=01   ,&   \sigma_{10}(B)&=10   ,\\
	\sigma_{01}(C)&=1    ,&   \sigma_{10}(C)&=1
.\end{align*}

In \cite{ABMP}, the following theorem is proved.
\begin{thrm}\label{thm:sigma-u}
An infinite ternary word $u\in\{A,B,C\}^\mathbb{N}$ is a 3iet word if and only if
 the words $\sigma_{01}(u)$ and $\sigma_{10}(u)$ are Sturmian.
\end{thrm}

This theorem motivated the authors of~\cite{AFMP} to introduce the relation of amicability of words.
\begin{dfntn}
	Let $w,w'\in\{0,1\}^*$, let $b\in\mathbb{N}$.
	We say that $w$ is \emph{$b$-amicable} to $w'$, if there exists a factor $v\in\{A,B,C\}^*$
	 of some 3iet word such that
	\[
		w=\sigma_{01}(v)
		,\qquad
		w'=\sigma_{10}(v)
		\quad\text{and}\quad
		\left|v\right|_B=b
	.\]
	We say that $w$ is \emph{amicable} to $w'$, if $w$ is $b$-amicable to $w'$ for some $b\in\mathbb{N}$,
	 and we denote it by $w\propto w'$.

	The ternary word $v$ is called a \emph{ternarization} of $w$ and $w'$,
	 and we write $v=\operatorname{ter}(w,w')$.
\end{dfntn}
It is easy to see that if $w\propto w'$, then they are factors of the same Sturmian word and
 their Parikh vectors coincide.

The ternarization is given uniquely for a pair $w$, $w'$.
For, let us see that if ternary words $v^{(1)}$, $v^{(2)}$ differ, then
 either $\sigma_{01}(v^{(1)})\neq\sigma_{01}(v^{(2)})$ or $\sigma_{10}(v^{(1)})\neq\sigma_{10}(v^{(2)})$.

In \cite{AFMP}, the notion of amicable words plays a crucial role in the enumeration of words with
 length $n$ occurring in a 3iet word.
In \cite{AMP-Mor}, the authors investigate ternary morphisms that have a non-degenerate 3iet fixed
 point using the following notion of amicability of two Sturmian morphisms.

\begin{dfntn}
	Let $\varphi,\psi$ be Sturmian morphisms over the alphabet $\{0,1\}$.
	We say that $\varphi$ is \emph{amicable} to $\psi$, if
	\begin{align*}
		\varphi(0)&\propto\psi(0)
		,\\
		\varphi(01)&\propto\psi(10)
		\\\text{and}\quad
		\varphi(1)&\propto\psi(1)
	.\end{align*}
	We denote this relation by $\varphi\propto\psi$.
	The morphism $\eta$ over the ternary alphabet $\{A,B,C\}$, given by
	\begin{align*}
		\eta(A)&=\operatorname{ter}\bigl(\varphi(0),\psi(0)\bigr)
		,\\
		\eta(B)&=\operatorname{ter}\bigl(\varphi(01),\psi(10)\bigr)
		,\\
		\eta(C)&=\operatorname{ter}\bigl(\varphi(1),\psi(1)\bigr)
	,\end{align*}
	 is called the \emph{ternarization} of morphisms $\varphi$ and $\psi$,
	 and is denoted by $\eta=\operatorname{ter}(\varphi,\psi)$.
	The set of these $\eta$ is denoted by $\mathcal{M}_{\mathrm{ter}}$.
\end{dfntn}

The ternarization of words is given uniquely by the words $u\propto v$, hence the ternarization
 of morphisms is given uniquely as well.

\begin{xmpl}
	Consider Sturmian morphisms $\varphi,\psi$ given by
	\begin{align*}
		\varphi(0)&=001
	,&
		\varphi(1)&=00101
	,&
		\psi(0)&=010
	,&
		\psi(1)&=01001
	.\end{align*}
	Then $\varphi\propto\psi$ and their ternarization $\eta=\operatorname{ter}(\varphi,\psi)$ satisfies
	\begin{align*}
		\eta(A)&=AB
	,&
		\eta(B)&=ABABB
	,&
		\eta(C)&=ABAC
	.\end{align*}
\end{xmpl}

The article \cite{AMP-Mor} states the following theorem:
\begin{thrm}\label{thm:eta-2}
	Let $\eta$ be a ternary morphism with non-degenerate 3iet fixed point.
	Then $\eta\in\mathcal{M}_{\mathrm{ter}}$ or $\eta^2\in\mathcal{M}_{\mathrm{ter}}$.
\end{thrm}

\section{Main results}

Analogously to the terminology introduced for Sturmian words and morphisms in \cite{BS94}, the ternarization $\eta$,
 having a 3iet fixed point, is \emph{locally 3iet-preserving}, i.e.\ there exists
 $u\in\mathcal{W}_{\mathrm{3iet}}$ such that $\eta(u)\in\mathcal{W}_{\mathrm{3iet}}$.
We now prove a partial result about \emph{(globally) 3iet-preserving} morphisms,
 i.e.\ ternary morphisms $\eta$ such that
\[
	\eta(u)\in\mathcal{W}_{\mathrm{3iet}} \quad\text{for all}\quad u\in\mathcal{W}_{\mathrm{3iet}}
.\]

\begin{prpstn}
	Let $\eta=\operatorname{ter}(\varphi,\psi)$ for amicable Sturmian morphisms $\varphi\propto\psi$.
	Then $\eta$ is a globally 3iet-preserving morphism.
\end{prpstn}

\begin{proof}
Directly from definitions we see that
\begin{align*}
	\sigma_{01}\eta(A)&=\varphi(0) ,& \sigma_{01}\eta(B)&=\varphi(01) ,& \sigma_{01}\eta(C)&=\varphi(1)
	,\\
	\sigma_{10}\eta(A)&=\psi(0) ,& \sigma_{10}\eta(B)&=\psi(10) ,& \sigma_{10}\eta(C)&=\psi(1)
.\end{align*}
Therefore
\begin{equation}\label{eq:sig-p}
	\sigma_{01}\eta(v)=\varphi\sigma_{01}(v)
	\qquad\text{and}\qquad
	\sigma_{10}\eta(v)=\psi\sigma_{10}(v)
\end{equation}
 for any factor $v$ of a 3iet word $u\in\mathcal{W}_{\mathrm{3iet}}$.
According to Theorem \ref{thm:sigma-u} we get that $\sigma_{01}(u)$ and $\sigma_{10}(u)$ are Sturmian words,
 and since $\varphi$ and $\psi$ are Sturmian morphisms, we obtain that $\sigma_{01}\eta(u)$ and
 $\sigma_{10}\eta(u)$ are Sturmian words as well, which means, according to the same theorem,
 that the word $\eta(u)$ is 3iet.
\end{proof}

\begin{prpstn}
	Let $\varphi_i\propto\psi_i$ be Sturmian morphisms, for $i=1,2$.
	Then
	\[
		\operatorname{ter}(\varphi_1,\psi_1)\circ\operatorname{ter}(\varphi_2,\psi_2)=\operatorname{ter}(\varphi_1\circ\varphi_2,\psi_1\circ\psi_2)
	.\]
\end{prpstn}

\begin{proof}
It can be shown that the relation of amicability is preserved by composition of morphisms.
More precisely $\varphi_1\varphi_2\propto\psi_1\psi_2$.
Denote $\eta_1=\operatorname{ter}(\varphi_1,\psi_1)$, $\eta_2=\operatorname{ter}(\varphi_2,\psi_2)$.
Using the relation \eqref{eq:sig-p}, we see that for all $v\in\{A,B,C\}^*$
\begin{align*}
	\sigma_{01}\eta_1\eta_2(v)&=\varphi_1\sigma_{01}\eta_2(v)=\varphi_1\varphi_2\sigma_{01}(v)
	\\\text{and}\quad
	\sigma_{10}\eta_1\eta_2(v)&=\psi_1\sigma_{10}\eta_2(v)=\psi_1\psi_2\sigma_{10}(v)
.\end{align*}
But this means that $\eta_1\eta_2=\operatorname{ter}(\varphi_1\varphi_2, \psi_1\psi_2)$.
\end{proof}

As a consequence of previous two propositions, we can state the following theorem.

\begin{thrm}
	The set $\mathcal{M}_{\mathrm{ter}}$ of all ternarizations of amicable Sturmian morphisms with the operation of composition of morphisms
	 is a sub-monoid of the monoid $\mathcal{M}_{\mathrm{3iet}}$ of all globally 3iet-preserving morphisms.
\end{thrm}

Unfortunately, $\mathcal{M}_{\mathrm{ter}}\subsetneqq\mathcal{M}_{\mathrm{3iet}}$. Consider for example the morphism
\begin{align}\label{eq:notter}
	\eta(A)&=B
	,&
	\eta(B)&=CAC
	,&
	\eta(C)&=C
.\end{align}
As shown in \cite{Hak}, this morphism is 3iet-preserving, but it can be easily verified that it is
 not a ternarization of any pair of Sturmian morphisms, using the following statement.

\begin{prpstn}
	A ternary morphism $\eta$ is a ternarization, i.e. $\eta\in\mathcal{M}_{\mathrm{ter}}$, if and only if it satisfies
	\[
		\sigma_{01}\eta(B)=\sigma_{01}\eta(AC)
		\quad\text{and}\quad
		\sigma_{10}\eta(B)=\sigma_{10}\eta(CA)
	.\]
\end{prpstn}

\begin{proof}
The implication $(\Rightarrow)$.
Suppose $\eta=\operatorname{ter}(\varphi,\psi)$.
According to \eqref{eq:sig-p} we get
\begin{gather*}
\sigma_{01}\eta(B)=\varphi\sigma_{01}(B)=\varphi(01)=\varphi\sigma_{01}(AC)=\sigma_{01}\eta(AC)
,\\
\sigma_{10}\eta(B)=\psi\sigma_{10}(B)=\psi(10)=\psi\sigma_{10}(CA)=\sigma_{10}\eta(CA)
.\end{gather*}

The implication $(\Leftarrow)$.
Define morphisms $\varphi$, $\psi$ as
\begin{align*}
	\varphi(0)&=\sigma_{01}\eta(A)
	,&
	\psi(0)&=\sigma_{10}\eta(A)
	,\\
	\varphi(1)&=\sigma_{01}\eta(C)
	,&
	\psi(1)&=\sigma_{10}\eta(C)
.\end{align*}
Immediately we get $\operatorname{ter}\bigl(\varphi(0),\psi(0)\bigr)=\eta(A)$ and $\operatorname{ter}\bigl(\varphi(1),\psi(1)\bigr)=\eta(C)$.
The words $\varphi(01)$ and $\psi(10)$ satisfy
\[
	\varphi(01)=\sigma_{01}\eta(AC)=\sigma_{01}\eta(B)
	\quad\text{and}\quad
	\psi(10)=\sigma_{10}\eta(CA)=\sigma_{10}\eta(B)
,\]
 which means that $\operatorname{ter}\bigl(\varphi(01),\psi(10)\bigr)=\eta(B)$.
\end{proof}

For the morphism \eqref{eq:notter}, we get $\sigma_{01}\eta(B)=010\neq011=\sigma_{01}\eta(AC)$.
Another even simpler example of a 3iet-preserving morphism that is not a ternarization
 is the morphism interchanging the letters $A$ and $C$.

Now, our goal will be to determine the number of amicable pairs of morphisms with incidence matrix $\mathbf A$ of $\det\mathbf A=\pm1$.
We will use the notion of $b$-amicable morphisms.
\begin{dfntn}
	Let $\varphi$ and $\psi$ be binary morphisms and let $b\in\mathbb{N}$.
	We say that $\varphi$ is \emph{$b$-amicable} to $\psi$,
	 if $\varphi$ is amicable to $\psi$ and the number of occurrences of $B$ in $\operatorname{ter}\bigl(\varphi(01),\psi(10)\bigr)$ is $b$.
\end{dfntn}

We now determine the numbers of pairs of $b$-amicable Sturmian morphisms.

\begin{prpstn}\label{proposition:cAb}
	Let $\mathbf A=\left(\begin{smallmatrix} p_0 & q_0 \\ p_1 & q_1 \end{smallmatrix}\right)\in\mathbb{N}^{2\times2}$
	 be a matrix with $\det\mathbf A=\pm1$ and $b\in\mathbb{N}$.
	Put $p=p_0+p_1$, $q=q_0+q_1$.
	Then the number $c_{\mathbf A}(b)$ of pairs of $b$-amicable morphisms with matrix $\mathbf A$ is equal to
	\begin{align*}
		c_{\mathbf A}(b)&=
		\begin{cases}
			\left\|\mathbf A\right\|-b & \text{if }\; \det\mathbf A=+1 \;\text{ and }\; 1\leq b\leq \min\{p,q\},
			\\
			\left\|\mathbf A\right\|-b-2 & \text{if }\; \det\mathbf A=-1 \;\text{ and }\; 0\leq b\leq \min\{p,q\}-1,
			\\
			0 & \text{otherwise,}
		\end{cases}
	\end{align*}
	where $\left\|\mathbf A\right\|=p+q$.
\end{prpstn}

First, let us state the following lemma.

\begin{lmm}\label{lemma:S}
	Let $\mathbf A=\left(\begin{smallmatrix} p_0 & q_0 \\ p_1 & q_1 \end{smallmatrix}\right)\in\mathbb{N}^{2\times2}$
	 be a matrix with $\det\mathbf A=\pm1$ and $b\in\mathbb{N}$.
	Put $p=p_0+p_1$, $q=q_0+q_1$ and $N=\left\|\mathbf A\right\|=p+q$.
	Let $S$ be a two-interval exchange with the slope $p/N$.
	Let $w^{(k)}$ be a word of the length $N$ that codes $S$ with the start point $k/N$,
	 for $k\in\{0,\dots,N-1\}$.
	
	Then $w^{(k)}$ is $b$-amicable to $w^{(\bar k)}$ if and only if $0\leq b\leq\min\{p,q\}$ and $\bar k-k=b$.	
\end{lmm}

\begin{proof}
Using \eqref{eq:Sxi}, we see that $S^i(k/N)\equiv (k-ip)/N\pmod1$, which is equivalent to
 $NS^i(k/N)\equiv k-ip\pmod N$.
We know that the numbers $p$ and $N$ are co-prime, thus
 the mapping $f_k:\{0,\dots,N-1\}\rightarrow\{0,\dots,N-1\}$ given by the congruence
 $f_k(i)\equiv k-ip\pmod N$ is a bijection.
As well, $f_{\bar k}(i)-f_k(i)\equiv \bar k-k\pmod N$.

Denote $m=\min\{p,q\}$ and $b=\bar k-k$.
Consider the following cases:
\begin{itemize}

\item Case $b<0$.
We shall see that $w^{(k)}$ is lexicographically larger than $w^{(\bar k)}$,
 i.e. if $i\in\mathbb{N}$ is the first position such that $w^{(k)}_i\neq w^{(\bar k)}_i$,
 then $w^{(k)}_i=1$ and $w^{(\bar k)}_i=0$.
Directly from the definition of amicability, if $w^{(k)}\propto w^{(\bar k)}$
 and $w^{(k)}\neq w^{(\bar k)}$, then $w^{(k)}$ is lexicographically smaller than $w^{(\bar k)}$.
These two facts make a contradiction.

\item Case $b\in\{0,\dots,m\}$.
Let $\mathcal{I}_a\subset\{0,\dots,N-1\}$ be a set of indices $i$ such that $w^{(k)}_i=a$ and $w^{(\bar k)}_i\neq a$,
 for both $a=0,1$.
To show that $w^{(k)}$ is $b$-amicable to $w^{(\bar k)}$, we need to show that
 $i\in \mathcal{I}_0$ implies $i+1\in \mathcal{I}_1$ and $\#\mathcal{I}_0=\#\mathcal{I}_1=b$.
The fact that $\left|w^{(k)}\right|_0=\left|w^{(\smash{\bar k})}\right|_0$ follows to $\#\mathcal{I}_0=\#\mathcal{I}_1$.

Let $i$ be an index such that $f_k(i)\in[p-b,p)$, thus $w^{(k)}_i=0$.
Then $f_{\bar k}(i)\in[p,p+b)$, thus $w^{(\bar k)}_i=1$.
This means $i\in \mathcal{I}_0$.
For these $i$, we have $f_k(i+1)\in[N-b,N)$ and $f_{\bar k}(i+1)\in[0,b)$, which means $i\in \mathcal{I}_1$.
There are exactly $b$ such indices $i$.

It remains to show that we covered the whole set $\mathcal{I}_0$.
Suppose $f_k(i)<p-b$, then $f_{\bar k}(i)<p$ and $w^{(\bar k)}_i=0$, which means $i\notin \mathcal{I}_0$.
Suppose $f_k(i)\geq p$, then $w^{(k)}_i=1$, which means $i\notin \mathcal{I}_0$.

\item Case $b\in\{m+1,\dots,N-m-1\}$.
Let $i$ be such index that $f_k(i)=p-1$.
Then $f_k(i+1)=N-1$.

If $p\leq q$, then $f_{\bar k}(i)=b+p-1$ and $f_{\bar k}(i+1)=b-1$,
 which means that $w^{(k)}_i w^{(k)}_{i+1}=01$ and $w^{(\bar k)}_i w^{(\bar k)}_{i+1}=11$.

If $p>q$, then $f_{\bar k}(i)=b-q-1$ and $f_{\bar k}(i+1)=b-1$,
 which means that $w^{(k)}_i w^{(k)}_{i+1}=01$ and $w^{(\bar k)}_i w^{(\bar k)}_{i+1}=00$.

Both these are in contradiction with $w^{(k)}\propto w^{(\bar k)}$.

\item Case $b\in\{N-m,\dots,N-1\}$.

Suppose $p<q$.
Then $j=2p$ solves the inequalities
\begin{align*}
	p&\leq j<N
	,&
	p&\leq j+b-N<N
	,\\
	p&\leq j-p<N
	,&
	0&\leq j+b-p-N<p
.\end{align*}
Let $i$ be an index such that $f_k(i)=j$.
Then the previous inequalities give $w^{(k)}_i w^{(k)}_{i+1}=11$ and $w^{(\bar k)}_i w^{(\bar k)}_{i+1}=10$,
 which is in a contradiction with $w^{(k)}\propto w^{(\bar k)}$.

\enlargethispage{12mm}

Suppose $p>q$.
Then $j=2p-b-1$ solves the inequalities
\begin{align*}
	0&\leq j<p
	,&
	0&\leq j+b-N<p
	,\\
	p&\leq j-p+N<N
	,&
	0&\leq j+b-p<p
.\end{align*}
Let $i$ be an index such that $f_k(i)=j$.
Then the previous inequalities give $w^{(k)}_i w^{(k)}_{i+1}=01$ and $w^{(\bar k)}_i w^{(\bar k)}_{i+1}=00$,
 which is a contradiction with $w^{(k)}\propto w^{(\bar k)}$.
 \qedhere

\end{itemize}
\end{proof}

\begin{proof}[Proof of Proposition \Rref{proposition:cAb}.]
Let $S$ be a 2-interval exchange transformation with the slope $\varepsilon=p/N$.
Let $k\in\mathbb{Z}$ and denote $w^{(k)}$ the word of the length $N=\left\|\mathbf A\right\|$ that codes the orbit
 of the point $\{k/N\}$ with respect to $S$.
From \cite{Seebold} we know that for every Sturmian morphism $\varphi$ with $\mathbf M_\varphi=\mathbf A$, there exists $k\in\{0,\dots,N-1\}$
 such that $\varphi(01)=w^{(k)}$, we will denote this morphism $\varphi^{(k)}$.

Let $\varphi_{\mathrm{std}}$ be a standard morphism with $\mathbf M_{\varphi_{\mathrm{std}}}=\mathbf A$.
Every Sturmian morphism $\varphi^{(k)}$ is a right conjugate to $\varphi_{\mathrm{std}}$,
 which means that there exist words $v,v'\in\{0,1\}*$ such that
\[
	\varphi^{(k)}(aa')=v01v'
	\quad\text{and}\quad
	\varphi^{(k)}(a'a)=v10v'
,\]
 where letters $a,a'$ satisfy $aa'=01$ for $\det\mathbf A=+1$ and $aa'=10$ for $\det\mathbf A=-1$.
This gives that $\varphi(aa')$ is 1-amicable to $\varphi(a'a)$.

Morphism $\varphi^{(k)}$ is $b$-amicable to $\varphi^{(\bar k)}$ if and only if the following conditions are satisfied:
\begin{enumerate}
\item \label{pf:cAb:1} $\varphi^{(k)}(01)$ is $b$-amicable to $\varphi^{(\bar k)}(10)$;
\item \label{pf:cAb:2} $\varphi^{(k)}(01)$ is amicable to $\varphi^{(\bar k)}(01)$;
\item \label{pf:cAb:3} Parikh vectors satisfy $\Psi(\varphi^{(k)}(0))=\Psi(\varphi^{(\bar k)}(0))$.
\end{enumerate}
The 2nd and 3rd conditions assures that $\varphi^{(k)}(0)\propto\varphi^{(\bar k)}(0)$ and $\varphi^{(k)}(1)\propto\varphi^{(\bar k)}(1)$.

Let us discuss the cases $\det\mathbf A=+1$ and $\det\mathbf A=-1$.
\begin{itemize}

\item Case $\det\mathbf A=+1$.
We know that $\varphi^{(k)}(01)$ is $1$-amicable to $\varphi^{(k)}(10)$,
 implying by Lemma \ref{lemma:S} that $\varphi^{(k)}(10)=w^{(k+1)}$.
This excludes $k=N-1$.

The 3rd condition is immediately satisfied by $\mathbf M_{\varphi^{(k)}}=\mathbf M_{\varphi^{(\bar k)}}$.
To satisfy the 1st condition, we need $(\bar k+1)-k=b$.
To satisfy the 2nd condition, we need $0\leq \bar k-k\leq\min\{p,q\}$.
These facts gives $0\leq k\leq \bar k\leq N-2$ and $1\leq b\leq\min\{p,q\}$,
 because the value $b=\min\{p,q\}+1$ is denied by Lemma \ref{lemma:S}.
For each admissible $b$, we have exactly $N-b$ pairs of indices $(k,\bar k)$.

\item Case $\det\mathbf A=-1$.
We know that $\varphi^{(k)}(10)$ is $1$-amicable to $\varphi^{(k)}(01)$,
 implying by Lemma \ref{lemma:S} that $\varphi^{(k)}(10)=w^{(k-1)}$.
This excludes $k=0$.

The 3rd condition is immediately satisfied by $\mathbf M_{\varphi^{(k)}}=\mathbf M_{\varphi^{(\bar k)}}$.
To satisfy the 1st condition, we need $(\bar k-1)-k=b$.
To satisfy the 2nd condition, we need $0\leq \bar k-k\leq\min\{p,q\}$.
These facts gives $1\leq k\leq \bar k\leq N-1$ and $0\leq b\leq\min\{p,q\}-1$,
 because the value $b=-1$ is denied by Lemma \ref{lemma:S}.
For each admissible $b$, we have exactly $N-b-2$ pairs of indices $(k,\bar k)$.
\qedhere

\end{itemize}
\end{proof}

\begin{rmrk}\label{rmrk:kbark}
The proof shows an interesting fact:
Suppose that
\begin{equation}\label{rmrk:kbark:eq:01}
\text{the word $\varphi^{(k)}(01)$ is $(b-\Delta)$-amicable to
$\varphi^{(\bar k)}(01)$}
\end{equation}
 and $c_{\mathbf A}(b)\neq0$.
Then the morphism $\varphi^{(k)}$ is $b$-amicable to $\varphi^{(\bar k)}$. The reason is as follows:
In the proof we considered all pairs of $(k, \bar k)$ and to satisfy 
\eqref{rmrk:kbark:eq:01} there is no other choice but $\bar k-k=b-\Delta$.
The condition $c_{\mathbf A}(b)\neq0$ is what we needed in the proof to show that $\varphi^{(k)}(01)$ is $b$-amicable to $\varphi^{(\bar k)}(10)$.
Thus the conditions 1, 2 from the proof are true; the condition 3 is straightforward.
\end{rmrk}

\begin{proof}[Proof of Theorem \Rref{thm:AmN}.]
The formula \eqref{eq:AmN} can be obtained by summation of numbers $c_{\mathbf A}(b)$
 from the previous proposition.
\end{proof}

To each pair of amicable Sturmian morphisms, an incidence matrix of its ternarization is assigned.
We now fully describe which matrices from $\mathbb{N}^{3\times3}$ are matrices of ternarizations.

\begin{thrm}\label{thm:B}
	A matrix $\mathbf B\in\mathbb{N}^{3\times3}$ is the incidence matrix of the ternarization of a pair of amicable Sturmian morphisms
	 if and only if there exists a matrix $\mathbf A=\left(\begin{smallmatrix} p_0 & q_0 \\ p_1 & q_1 \end{smallmatrix}\right)\in\mathbb{N}^{2\times2}$
	  with $\det\mathbf A=\Delta=\pm1$ and numbers $b_0,b_1\in\mathbb{N}$ such that
	\begin{itemize}
	\item[(a)]
		$\left|\frac{b_0(p_1+q_1)-b_1(p_0+q_0)}{p_0+q_0+p_1+q_1}\right|<1$,
	\item[(b)]
		$\frac{1-\Delta}2\leq b_0+b_1\leq\min\{p_0+p_1,q_0+q_1\}-\frac{\Delta+1}2$,
	\item[(c)]
		$\mathbf B=\mathbf P
		\left(\!\begin{smallmatrix}
			{\textstyle\mathbf A} & \begin{smallmatrix}b_0\\b_1\end{smallmatrix} \\ 0\;\,0 & \Delta
		\end{smallmatrix}\!\right)
		\mathbf P^{-1}$,
		where
		$\mathbf P=\left(\!\begin{smallmatrix}
			1&0&0 \\ 1&1&1 \\ 0&1&0
		\end{smallmatrix}\!\right)$.
	\end{itemize}
\end{thrm}

\enlargethispage{8mm}

\begin{proof}[Proof of the implication $(\Rightarrow)$.]
Let us denote $p=p_0+p_1$, $q=q_0+q_1$, $N=p+q$ and $b=b_0+b_1+\Delta$.
Then we can see that condition (c) gives
\begin{equation}\label{eq:B}
	\mathbf B=\begin{pmatrix}p_0-b_0 & b_0 & q_0-b_0 \\ p-b & b & q-b \\ p_1-b_1 & b_1 & q_1-b_1\end{pmatrix}
.\end{equation}

The fact that (c) is necessary for $\mathbf B$ to be an incidence matrix of a ternarization is shown in \cite[Remark 13]{AMP-Mat}.
Condition (b) is necessary according to Proposition \ref{proposition:cAb},
 so we only need to show that (a) is satisfied for the matrix of the ternarization $\eta=\operatorname{ter}(\varphi,\psi)$
 of a pair of amicable Sturmian morphisms $\varphi\propto\psi$.

We can see that $\mathbf A=\left(\begin{smallmatrix} p_0 & q_0 \\ p_1 & q_1 \end{smallmatrix}\right)$
 is necessarily an incidence matrix of both $\varphi$ and $\psi$.
Let $S$ be a 2-interval exchange transformation with a rational slope $\varepsilon=p/N$.
Then there exist numbers $k, \bar k\in\{0,\ldots,N-2\}$ such that $\varphi(01)$, $\psi(01)$ code transformation $S$
 with start points $x_0=k/N$, $\bar x_0=\bar k/N$, respectively; moreover, $\bar k-k=b-\Delta$.
We need to determine the value of $b_0=\left|\operatorname{ter}\bigl(\varphi(0),\psi(0)\bigr)\right|_B$.
The number $b_0$ is equal to the number of indices $i\in\{0,1,\ldots,p_0+q_0-1\}$ such that $S^ix_0\in\bigl[(p-b+\Delta)/N,p/N\bigr)$,
 because for exactly these $i$, we have $S^ix_0<p/N\leq S^i\bar x_0$.

Let $X=\bigl\{\{x_0-ip/N\}\big|i\in\mathbb{N}, 0\leq i<p_0+q_0\bigr\}$.
Put $p'=p+\Delta/(p_0+q_0)$, and
let $Y=\bigl\{\{x_0-ip'/N\}\big|i\in\mathbb{N}, 0\leq i<p_0+q_0\bigr\}$.
We can see that $0\leq\Delta\bigl((x_0-ip/N)-(x_0-ip'/N)\bigr)=i/(p_0+q_0)N<1/N$.
Thus $x_0-ip/N\in \bigl[\tfrac{p-b+\Delta}N,\tfrac pN\bigr)$ if and only if
\begin{equation}\label{eq:xipN}
	x_0-ip'/N\in\begin{cases}
 	\bigl(\tfrac{p-b}N,\tfrac{p-1}N\bigr] &\text{in the case }\; \Delta=+1
 	,\\
 	\bigl[\tfrac{p-b-1}N,\tfrac pN\bigr) &\text{in the case }\; \Delta=-1
 	.\end{cases}
\end{equation}
In both cases, the length of the interval is $\tfrac{b-\Delta}N$.
From $\Delta=\det\mathbf A=\det\left(\begin{smallmatrix} p_0 & p_0+q_0 \\ p & N \end{smallmatrix}\right)$,
 it is easy to see that
\[
	\frac{p'}{N}
	=\frac{p+\Delta/(p_0+q_0)}{N}
	=\frac pN + \frac{p_0N-p(p_0+q_0)}{N(p_0+q_0)}
	=\frac{p_0}{p_0+q_0}
.\]
Because $p_0$ is co-prime to $p_0+q_0$, we get
 $\bigl\{\{i p_0/(p_0+q_0)\}\big|i\in\mathbb{N}, 0\leq i<p_0+q_0\bigr\}=\bigl\{i/(p_0+q_0)\big|i\in\mathbb{N}, 0\leq i<p_0+q_0\bigr\}$.
But this means that the set $Y$ is uniformly distributed on the interval $[0,1)$, therefore
\[
	b_0=\#\Bigl( X \cap \bigl[\tfrac{p-b+\Delta}N,\tfrac pN\bigr) \Bigl)
	\in\bigl\{\lfloor\beta\rfloor,\lceil\beta\rceil\bigr\}
,\]
where $\beta=(p_0+q_0)\frac{b-\Delta}N$ is number of elements of $Y$ multiplied by
 the length of the interval \eqref{eq:xipN}.
Together we get
\begin{equation}\label{eq:bb}
	\left|\beta-b_0\right|<1
,\end{equation}
which is equivalent to condition (a).
\end{proof}

The proof of the other implication is divided into several lemmas.

\enlargethispage{-5mm}

\begin{lmm}\label{lemma:1}
Let $\mathbf A=\left(\begin{smallmatrix} p_0 & q_0 \\ p_1 & q_1 \end{smallmatrix}\right)\in\mathbb{N}^{2\times2}$ with $\det\mathbf A=\Delta=\pm1$,
 let $b\in\mathbb{N}$ with $\frac{1+\Delta}2\leq b\leq\min\{p_0+p_1,q_0+q_1\}-\frac{1-\Delta}2$.

Denote $N=\left\|\mathbf A\right\|$, $p=p_0+p_1$ and $q=q_0+q_1$ integers, $I=\bigl[\frac{p-b+\Delta}N,\frac pN\bigr)$ an interval,
 $X_k=\bigl\{\{k/N\}, S\{k/N\}, S^2\{k/N\}, \ldots, S^{p_0+q_0-1}\{k/N\}\bigr\}$ a set of numbers for any $k\in\mathbb{Z}$,
 where $S$ is the 2-interval exchange with the slope $\varepsilon=p/N$,
 and denote $\beta=\frac{p_0+q_0}N(b-\Delta)$.

Then for all $b_0\in\bigl\{\lfloor\beta\rfloor,\lceil\beta\rceil\bigr\}$ such that 
\begin{equation}\label{eq:1bpq}
	b_0\leq\min\{p_0,q_0\}
	\quad\text{and}\quad
	b-\Delta-b_0\leq\min\{p_1,q_1\}
,\end{equation}
there exist $k',k''\in\{0,\ldots,N-1\}$, $k'\neq k''$ such that
\begin{equation}\label{eq:k''}
	\#(X_{k'}\cap I)=\#(X_{k''}\cap I)=b_0
.\end{equation}
\end{lmm}

\begin{proof}
Denote $r(k)=\#(X_{k}\cap I)$ for $k\in\mathbb{Z}$.
We can see that $\sum_{k=0}^{N-1} r(k)=(b-\Delta)(p_0+q_0)$.
According to \eqref{eq:bb}, we know that
 $r(k)\in\bigl\{\lfloor\beta\rfloor,\lceil\beta\rceil\bigr\}$ for all $k\in\mathbb{Z}$.
Let
\begin{align*}
	C_L&=\#\bigl\{ k\in\{0,\ldots,N-1\} \big| r(k)=\lfloor\beta\rfloor \bigr\}
	,\\
	C_U&=\#\bigl\{ k\in\{0,\ldots,N-1\} \big| r(k)=\lceil\beta\rceil \bigr\}
.\end{align*}
These numbers satisfy the equations
\begin{align}\notag
	C_L\lfloor\beta\rfloor+C_U\lceil\beta\rceil&=N\beta
	\\\text{and}\qquad\label{l:3:eq:LUN}
	C_L+C_U&=N
.\end{align}

If $C_L=0$ or $C_U=0$, necessarily $\beta\in\mathbb{N}$ and \eqref{eq:k''} is satisfied for all $k\in\mathbb{Z}$.

If $C_L\geq2$, we have two different $k\in\mathbb{Z}$ satisfying \eqref{eq:k''} for $b_0=\lfloor\beta\rfloor$.
Similarly if $C_U\geq2$, we have two different $k\in\mathbb{Z}$ satisfying \eqref{eq:k''} for $b_0=\lceil\beta\rceil$.

We will show that $C_L=1$ implies $\lfloor\beta\rfloor$ not to satisfy the condition \eqref{eq:1bpq},
and similarly for $C_U$ and $\lceil\beta\rceil$.

If $C_U$ and $C_L$ are non-zero then there is a unique solution 
\[
	C_L=N\{-\beta\}
	\quad\text{and}\quad
	C_U=N\{\beta\}
.\]
Using relation $p_0N-(p_0+p_1)(p_0+q_0)=\Delta$, we get
\begin{align}\notag
	C_U&\equiv(p_0+q_0)(b-\Delta) \pmod N
	\\\label{l:3:eq:bD}
	b-\Delta&\equiv-\Delta (p_0+p_1)C_U \pmod N
.\end{align}

Let us suppose $C_U=1$ or $C_L=1$, i.e. $C_U\equiv\pm1\pmod N$ due to \eqref{l:3:eq:LUN}.
Then \eqref{l:3:eq:LUN} and \eqref{l:3:eq:bD} lead to $b=(p_0+p_1)+\Delta$ or $b=(q_0+q_1)+\Delta$.
For $\Delta=+1$, this is in contradiction with the conditions.
For $\Delta=-1$, discuss the following two cases.
\begin{itemize}

\item
Case $b=(p_0+p_1)+\Delta$.
This happens when $C_U=1$.
But it means that $b_0=\lceil\beta\rceil$ is equal to $\bigl\lceil\frac{p_0N-\Delta}N\bigr\rceil=p_0+1$
 and this case is excluded by the condition \eqref{eq:1bpq}.

\item
Case $b=(q_0+q_1)+\Delta$.
This happens when $C_L=1$.
But it means that $b_0=\lfloor\beta\rfloor$ is equal to $q_0-1$ hence $b-\Delta-b_0=q_1+1$,
 which is excluded by \eqref{eq:1bpq}.
\qedhere

\end{itemize}
\end{proof}

\begin{lmm}\label{lemma:2}
Let us have the same hypothesis as in Lemma \ref{lemma:1}.

Define morphisms $\varphi_k$ for $k\in\mathbb{Z}$ in the following way:
\begin{itemize}
\item the word $\varphi_k(0)$ codes $\{k/N\}, S\{k/N\}, \ldots, S^{p_0+q_0-1}\{k/N\}$;
\item the word $\varphi_k(1)$ codes $S^{p_0+q_0}\{k/N\}, \ldots, S^{N-1}\{k/N\}$.
\end{itemize}
Let $k_0\in\mathbb{Z}$ be such integer that $\#(X_{k_0}\cap I)=\#(X_{k_0-p}\cap I)$.
Then
\[
	\varphi_{k_0}\propto\varphi_{k_0+b-\Delta}
	\quad\text{or}\quad
	\varphi_{k_0-p}\propto\varphi_{k_0-p+b-\Delta}
,\]
and the number of B's in the ternarization of the images of the letter 0 is $\#(X_{k_0}\cap I)$.
\end{lmm}

\begin{proof}
Let $k\in\mathbb{Z}$ and let us consider the orbit 
\begin{equation}\label{eq:S-1}
	\{k/N\}, S\{k/N\}, \ldots, S^{p_0+q_0-1}\{k/N\}
.\end{equation}
Let $t^{(k)}$ be a word of the length $p_0+q_0$ that codes \eqref{eq:S-1}
 to the alphabet $\{0,0',1,1'\}$ with the following code:
\begin{equation}\label{eq:tki}
	t^{(k)}_i=\begin{cases}
	0 &\text{if }\; S^i\{k/N\}\in\bigl[0,\frac{p-b+\Delta}N\bigl)
	,\\
	0' &\text{if }\; S^i\{k/N\}\in\bigl[\frac{p-b+\Delta}N,\frac pN\bigr)=I
	,\\
	1 &\text{if }\; S^i\{k/N\}\in\bigl[\frac pN,\frac{N-b+\Delta}N\bigr)
	,\\
	1' &\text{if }\; S^i\{k/N\}\in\bigl[\frac{N-b+\Delta}N,1\bigr)
	. \end{cases}
\end{equation}
From definition of $S$, we see that $t^{(k)}_i=0'\Leftrightarrow t^{(k)}_{i+1}=1'$.
Define two morphisms $\tau, \tau':\{0,0',1,1'\}^*\rightarrow\{0,1\}^*$ as
\begin{align*}
	\tau(0)&=0 ,& \tau(0')&=0 ,& \tau(1)&=1 ,& \tau(1')&=1
	,\\
	\tau'(0)&=0 ,& \tau'(0')&=1 ,& \tau'(1)&=1 ,& \tau(1')&=0
.\end{align*}


If $t^{(k)}$ does not start with $1'$ and does not end with $0'$,
 then the word $\varphi_k(0)=\tau(t^{(k)})$ is
 $\left|t^{(k)}\right|_{0'}$-amicable to $\tau'(t^{(k)})=\varphi_{k+b-\Delta}(0)$.
Moreover, $\left|t^{(k)}\right|_{0'}=\#(X_{k}\cap I)$.
To show this, notice that $S\{k_0/N\}=\{(k_0-p)/N\}$, which means that there exist letters $a,a'\in\{0,0',1,1'\}$
 such that $t^{(k_0)}a=a't^{(k_0-p)}$ and $a=0'\Leftrightarrow a'=0'$,
 because the numbers of letters $0'$ in the words $t^{(k_0)}$ and $t^{(k_0-p)}$ coincide.

Consider these two cases:
\begin{itemize}
\item If $a=0'$ then the last letter of $t^{(k_0)}$ is not $0'$ since this implies $a'=1'$.
 This~yields $\varphi_k(0)\propto\varphi_{k+b-\Delta}(0)$ for $k=k_0$.
\item If $a\neq0'$ then $t^{(k_0-p)}$ does not start with $1'$ and does not end with $0'$.
 This~yields $\varphi_k(0)\propto\varphi_{k+b-\Delta}(0)$ for $k=k_0-p$.
\end{itemize}

Similar reasoning leads to the amicability of the images of the letter $1$.
Thus by concatenation $\varphi_k(01)\propto\varphi_{k+b-\Delta}(01)$.
The condition on $b$ is the same as in Proposition
 \ref{proposition:cAb}, hence Remark \ref{rmrk:kbark} applies.
\end{proof}

\begin{lmm}\label{lemma:3}
Let us have the same hypothesis as in Lemma \ref{lemma:1}.

Let $k_0\in\mathbb{Z}$ be a number such that if $\Delta=-1$ and $b=\min\{p,q\}-1$ then
\begin{equation}\label{eq:k0ne}
	k_0\not\equiv\begin{cases}
	\phantom{p-b}-1\pmod N\quad\text{in the case $p>q$}
	,\\
	p-b-1\pmod N\quad\text{in the case $p<q$}
	.\end{cases}
\end{equation}
Then
\[
	\#(X_{k_0}\cap I)=\#(X_{k_0+p}\cap I)
	\quad\text{or}\quad
	\#(X_{k_0}\cap I)=\#(X_{k_0-p}\cap I)
.\]
\end{lmm}

\enlargethispage{10mm}

\begin{proof}
Define the words $t^{(k)}$ by \eqref{eq:tki} in the same way as in the previous proof.
Denote $\ell=p_0+q_0$.
Then we know that there exist letters $a_0, \ldots, a_{\ell+1}\in\{0,0',1,1'\}$ such that
\begin{align*}
	t^{(k_0+p)}&=a_0a_1a_2\cdots a_{\ell-1}
	,\\
	t^{(k_0)}&=\phantom{a_0}a_1a_2\cdots a_{\ell-1}a_\ell
	,\\
	t^{(k_0-p)}&=\phantom{a_0a_1}a_2\cdots a_{\ell-1}a_\ell a_{\ell+1}
.\end{align*}
Let us remind that $\#(X_k\cap I)=\left|t^{(k)}\right|_{0'}$.
The proof will be done by contradiction.
Suppose that $\left|t^{(k_0+p)}\right|_{0'}\neq\left|t^{(k_0)}\right|_{0'}\neq\left|t^{(k_0-p)}\right|_{0'}$.
There are only two possible values of these numbers, thus $\left|t^{(k_0+p)}\right|_{0'}=\left|t^{(k_0-p)}\right|_{0'}$.
This together gives either $a_0=a_{\ell+1}=0'$ or $a_1=a_\ell=0'$.
It means that there exist $\xi\in I=\bigl[\frac{p-b+\Delta}N,\frac pN\bigr)$ and $\omega\in\{+1,-1\}$
 such that $S^{\ell+\omega}\xi\in I$.
Without the loss of generality $\xi\in\frac 1N\mathbb{Z}$.
Since $\ell p=p_0N-\Delta$, we have
\[
	S^{\ell+\omega}\xi
	\equiv \xi-\frac{(\ell+\omega)p}N
	\equiv \xi+\frac{\Delta-\omega p}N
	\pmod 1
.\]
Because $\left|S^{\ell+\omega}\xi-\xi\right|<1$ we have
\begin{align*}
	& S^{\ell+\omega}\xi-\xi=\frac{\Delta-\omega p}N
	\\\text{or}\quad
	& S^{\ell+\omega}\xi-\xi=\frac{\Delta-\omega p}N+\omega=\frac{\Delta+\omega q}N
,\end{align*}
since $1-p/N=q/N$.
This enforces $b-1-\Delta\geq\min\{p,q\}-1$ for the interval $I$ to be large enough
 to contain both $\xi$ and $S^{\ell+\omega}\xi$.

For $\Delta=+1$, this is in contradiction with $b\leq\min\{p,q\}$.

For $\Delta=-1$ we get only one admissible $b=\min\{p,q\}-1$.
The case $p=\min\{p,q\}$ means $\omega=-1$ and $\xi=\frac{p-b-1}N$, which implies $k_0\equiv p-b-1\pmod N$.
The case $q=\min\{p,q\}$ means $\omega=+1$ and $\xi=\frac{p-1}N$, which implies $k_0\equiv -1\pmod N$.
Both these cases are excluded by \eqref{eq:k0ne}.
\end{proof}

\begin{proof}[Proof of the implication $(\Leftarrow)$.]
From \cite[Remark 13]{AMP-Mat}, the incidence matrix of the ternarization $\operatorname{ter}(\varphi,\psi)$ is fully described by
 the matrix $\mathbf A$ and numbers $b_0$ and $b=b_0+b_1+\Delta$.
The condition (a) is equivalent to \eqref{eq:bb} and it gives at most two values of $b_0$.
If $\beta\in\mathbb{N}$, there is nothing to do as we have at least one pair of $b$-amicable morphisms
 $\varphi\propto\psi$ for $\mathbf A$, and its incidence matrix satisfies all three conditions.

For $\beta\notin\mathbb{N}$, we want to show that for both $b_0\in\bigl\{\lfloor\beta\rfloor,\lceil\beta\rceil\bigr\}$
 there exist $\varphi\propto\psi$ with $\left|\operatorname{ter}\bigl(\varphi(0),\psi(0)\bigr)\right|_B=b_0$.
Because the elements of the matrix $\mathbf B$ are non-negative, the condition \eqref{eq:1bpq} of Lemma \ref{lemma:1}
 is satisfied and we have two different $k', k''$.
At least one of them satisfies \eqref{eq:k0ne}.
Lemma \ref{lemma:3} then provides $k_0$ satisfying the conditions of Lemma \ref{lemma:2} that gives
 a pair of amicable Sturmian morphisms, ternarization of which has the incidence matrix $\mathbf B$.
\end{proof}

\section{Conclusions and open problems}

Matrices of 3iet-preserving morphisms were studied in \cite{AMP-Mat}.
The authors give a necessary condition on $\mathbf B\in\mathbb{N}^{3\times3}$ to be an incidence matrix
 of a 3iet-preserving morphism:
\[
            \mathbf B \mathbf E \mathbf B^{\mathsf T} = \pm \mathbf E,
            \quad\text{where}\quad
            \mathbf E=\left(\begin{smallmatrix}
            	\phantom{-}0 & \phantom{-}1 & \phantom{-}1
            	\\
            	-1 & \phantom{-}0 & \phantom{-}1
            	\\
            	-1 & -1 & \phantom{-}0
            \end{smallmatrix}\right)
.\]
However, this condition is not sufficient.
In our contribution, we study 3iet-preserving morphisms $\eta=\operatorname{ter}(\varphi,\psi)$
 arising from pairs of amicable Sturmian morphisms $\varphi\propto\psi$.
Our Theorem \ref{thm:B} gives sufficient and necessary condition
 for any matrix $\mathbf B\in\mathbb{N}^{3\times3}$ to satisfy $\mathbf B=\mathbf M_\eta$ for some ternarization $\eta=\operatorname{ter}(\varphi,\psi)$.

It remains to answer the question about the role of the monoid
\[
            \mathcal{M}_{\mathrm{ter}}=\bigl\{\operatorname{ter}(\varphi,\psi)\big|\varphi,\psi \text{ amicable morphisms}\bigr\}
\]
 in the whole monoid $\mathcal{M}_{\mathrm{3iet}}$ of all 3iet-preserving morphisms.
It seems that using similar proof as for Theorem \ref{thm:eta-2} (see \cite{AMP-Mor}) we can prove
 the following statement.
\begin{cnjctr}
Let $\eta\in\mathcal{M}_{\mathrm{3iet}}$.
Then one of $\eta$, $\eta\circ\xi_1$, $\eta\circ\xi_2$ or $\eta\circ\xi_1\circ\xi_2$ is in $\mathcal{M}_{\mathrm{ter}}$, where
\begin{align*}
	\xi_1(A)&=C ,& \xi_1(B)&=B ,& \xi_1(C)&=A ,
\\
	\xi_2(A)&=B ,& \xi_2(B)&=ACA ,& \xi_2(C)&=A .
\end{align*}
\end{cnjctr}



\end{document}